\documentclass[12pt,leqno]{article}
\usepackage[pdftex]{graphicx}
\usepackage{amsfonts, color, amsmath}

\newcounter{conjecture}
\newcounter{remark}
\newcounter{corollary}

\newtheorem{corollary}{Corollary}
\newtheorem{theorem}{Theorem}
\newtheorem{lemma}{Lemma}
\newtheorem{proposition}{Proposition}

\newcommand{\und}{\underline}

\newcommand{\lar}{\longrightarrow}
\newcommand{\eps}{\varepsilon}

\newcommand{\aaa}{\alpha}

\newcommand{\lll}{\label}
\newcommand {\rrr}[1]{(\ref{#1})}

\def \pf{{\bf Proof: }}

\def \be{\begin{equation}}
\def \ee{\end{equation}}
\def \bt{\begin{theorem}}
\def \et{\end{theorem}}
\def \bc{\begin{corollary}}
\def \ec{\end{corollary}}
\def \bea{\begin{eqnarray}}
\def \eea{\end{eqnarray}}
\def \bas{\begin{eqnarray*}}
\def \eas{\end{eqnarray*}}

\def \noi{\noindent}
\def \aa{\alpha}

\def \DD{\Delta}

\def \la{\lambda}

\def \th{\theta}

\def \pp{\phi}

\def \bbs{\backslash}
\def \vski{\vspace{12pt}}
\def \ff{\infty}

\def \CC{{\mathbb C}}
\def \DDD{{\mathbb D}}

\def \GG{{\cal G}}

\def \RR{{\mathbb R}}

\def \XX{{\cal X}}

\def \({\left(}
\def \){\right)}

\def \nn{\nonumber}

\def \bc{\begin{center} }
\def \ec{\end{center} }
\def \bs{\begin{slide} }
\def \es{\end{slide} }

\def\square{{\vcenter{\vbox{\hrule height.3pt
         \hbox{\vrule width.3pt height5pt \kern5pt
            \vrule width.3pt}
         \hrule height.3pt}}}}
\def\qed{{\hfill $\square$ \bigskip}}

\newcounter{cccases}
\begin{document}

\title{Simple random walk on distance-regular graphs}

\author{Greg Markowsky}

\maketitle

\begin{abstract}
A survey is presented of known results concerning simple random walk on the class of distance-regular graphs. One of the highlights is that electric resistance and hitting times between points can be explicitly calculated and given strong bounds for, which leads in turn to bounds on cover times, mixing times, etc. Also discussed are harmonic functions, moments of hitting and cover times, the Green's function, and the cutoff phenomenon. The main goal of the paper is to present these graphs as a natural setting in which to study simple random walk, and to stimulate further research in the field.
\end{abstract}

\tableofcontents

\section{Introduction}

The book \cite{doysne} documents a connection between electric resistance in circuits and simple random walk on graphs. There a large number of elegant results are proved, but the results in many cases are tempered by the difficulty of calculating explicit resistances for specific graphs, even those with relatively few edges. In essence, resistance may in many cases allow us to reformulate a difficult problem as a different but equally difficult one. It is natural, therefore, to concentrate on a class of graphs for which we may perform calculations easily. The {\it distance-regular graphs} is a large class of such graphs, and these form the topic of this paper. The definition of distance-regular graphs is restrictive enough so that such graphs are richly structured yet permissive enough to allow for a wealth of examples. The structure of these graphs allow for a number of concise results on resistance to be deduced, and these are the starting point of our investigation. We will see that these results lead to good estimates for such quantities as hitting times, cover times, and mixing times for distance-regular graphs.

\vski

This paper is intended as a survey, albeit a survey of a field that may well be in a relatively immature state. Nevertheless, an effort has been made to include all known general results concerning random walks on distance-regular graphs. The paper does prove a number of new results; however, the proofs contained herein are all relatively simple. The more difficult results, such as Theorems \ref{bigguy}, \ref{p1}, \ref{sympdev}, \ref{gcrazy}, \ref{thor}, \ref{cleo}, \ref{varad}, \ref{clint}, and \ref{bels}, as well as results connecting resistance to random walk, appear in other places.  It is hoped that other researchers will be suitably interested to work in this field themselves.



\subsection{Distance-regular graphs}

All the graphs considered in this paper are finite, undirected and
simple (for unexplained terminology and more details, see for example \cite{drgraphs}). We begin by fixing notation. Let $G$ be a
connected graph; technically $G$ is a vertex set together with an edge set, but for ease of notation in this paper we will associate $G$ with its set of vertices. We will write $x \sim y$ if there is an edge connecting $x$ to $y$. If there is an integer $k$ such that each vertex of a graph has exactly $k$ neighbors we say that $G$ is {\it regular} of {\it degree k}. The distance $d(x,y)$ between
any two vertices $x,y$ of $G$
is the length of a shortest path between $x$ and $y$ in $G$. An {\it automorphism} of $G$ is a bijection $\th$ from $G$ to itself such that $x \sim y$ if and only if $\th(x) \sim \th(y)$ for all $x,y \in G$. A {\it distance-transitive graph} is a graph for which if we are given two pairs of vertices $x_1,y_1$ and $x_2,y_2$ with $d(x_1,y_1) = d(x_2,y_2)$ then we can always find an automorphism $\th$ of the graph such that $\th(x_1)=x_2$ and $\th(y_1)=y_2$. It is clear that nontrivial distance-transitive graphs are highly structured. Distance-regular graphs are natural combinatorial generalizations of distance-transitive graphs, as we now describe.

\vski

For a vertex $x \in G$, define $K_i(x)$ to be the set of
vertices which are at distance $i$ from $x~(0\le i\le
D)$ where $D:=\max\{d(x,y)\mid x,y\in G\}$ is the diameter
of $G$. In addition, define $K_{-1}(x):=\emptyset$ and $K_{D+1}(x)
:= \emptyset$. A connected graph $G$ with diameter $D$ is called
{\em distance-regular} if there are integers $b_i,c_i$ $(0 \le i
\leq D)$ such that for any two vertices $x,y \in G$ with $d(x,y)=i$, there are precisely $c_i$
neighbors of $y$ in
$K_{i-1}(x)$ and $b_i$ neighbors of $y$ in $K_{i+1}(x)$
(cf. \cite[p.126]{drgraphs}). It is immediate that a distance-regular graph $G$ is regular with degree
$k := b_0$, and we define $a_i:=k-b_i-c_i$ for notational convenience. It is clear that the class of distance-transitive graphs is contained in the class of distance-regular graphs, and in fact this containment is proper. We refer to the list $(b_0,b_1,\ldots ,b_{D-1};c_1,c_2,\ldots,c_D)$ as the {\it intersection array} of $G$. Note that in all cases $c_0=b_D=0$, and $a_i=b_0-b_i-c_i$, so the intersection array implicitly contains the required information on the $a_i$'s as well. If we fix a vertex $x$ of $G$, then $|K_i(x)|$ does not depend on the
choice of $x$ as $c_{i+1} |K_{i+1}(x)| =
b_i |K_i(x)|$ holds for $i =1, 2, \ldots D-1$, so that in fact $|K_i(x)| = \frac{b_0 \ldots b_{i-1}}{c_1 \ldots c_i}$.  



\vski

This simple definition forces a great deal of structure upon any distance-regular graph. That these form a natural class upon which to study random walk is indicated by the fact that, if $X_t$ is a simple random walk on a distance-regular graph $G$ and $u$ is a prescribed vertex in $G$, then $Y_t = d(X_t,u)$ is itself a Markov chain; this idea is studied in detail in Section \ref{proj}. We will see in fact that the intersection array contains the information required to calculate all electric resistances between points, and therefore captures many important quantities related to simple random walk.
The study of these graphs is a large field in itself, and many researchers have worked to extend the known library and to characterize the intersection arrays which occur. \cite{drgraphs} contains numerous examples, but we will give a small number in the next section.

\vski

In all ensuing sections, except where otherwise stated, any graph referenced will be a distance-regular graph.

\subsection{Historical overview and examples} \lll{hist}

Distance-regular graphs were defined by Biggs in the late 1960's as a combinatorial generalization of distance-transitive graphs. In the early 1970's Delsarte in his thesis \cite{Delsarte} introduced association schemes in order to study codes by algebraic methods. It turned out that the so called $P$-polynomial association schemes as introduced by Delsarte are exactly the same objects as the distance-regular graphs introduced by Biggs.  Most of the known examples of large diameter come from classical constructions, and are highly symmetric. A large number of these graphs appear naturally in theoretical computer science and coding theory. Here are a few simple examples which it is hoped will be of interest to readers without prior knowledge of the field.
\\

(i) The {\it complete graph} $K_n$, which is the graph with $n$ vertices and each vertex adjacent to all others. These have intersection array $(n-1;1)$. Many of the propositions of this paper reduce to trivialities in the case of complete graphs. \\

(ii) The {\it cycle graphs} $C_n$, which are the unique regular graphs of degree 2 with $n$ vertices. These have intersection arrays $(2,1, \ldots ,1; 1, \ldots, 1, 2-\{n(\mbox{mod }2)\})$. Simple random walks on cycles present points of interest (see for instance \cite[p. 25]{pers}); however, the primary methods and results developed in this paper do not have much application to this case--see the remark following Proposition \ref{bigguy} below--and we will for this reason exclude them from much of our analysis. \\

(iii) The {\it Hamming graphs}, $H(m,q)$. Let $Q$ be a set of size $q$. The vertex set of $H(m,q)$ are the vectors of length $m$ with entries in $Q$ and two such vectors are adjacent if they differ in exactly one position. The graph $H(m,2)$ is often referred to as the {\it hypercube graph} $Q_m$. $H(m,q)$ has intersection array $(m(q-1),(m-1)(q-1),(m-2)(q-1), \ldots , q-1; 1,2,3, \ldots , m)$. The Hamming graphs have been found to have interesting connections with the study of error-correcting codes. \\

(iv) The {\it Johnson Graphs} $J(m,q)$. Let $N$ be a set of size $m$ and let the vertices of $J(m,q)$ be the set of subsets of $N$ of size $q$. Two such subsets are adjacent if their intersection has cardinality $q-1$. The graph $J(m,q)$ has intersection array $((m-q)q,(m-q-1)(q-1), \ldots; 1,4,9, \ldots)$. The Johnson graphs arise naturally in the study of the {\it Johnson scheme}, which became well-known in connection with coding theory. \\

(v) The {\it Grassman Graphs} $J_q(m,t)$. Let $V$ be a vector space of dimension $m$ over the finite field $F$ with $q$ elements ($q$ a prime power). The vertices are the subspaces of dimension $t$ (over $F$) of $V$ and such two vertices are adjacent if their intersection is a vector space of dimension $t-1$. Grassman graphs have been found to be of interest in relation to quantum physics. \\

(vi) The {\it Odd graphs} $O_m$. The vertex set may be taken to be all subsets of size $m-1$ of a set of size $2m-1$, with two vertices being adjacent precisely when the corresponding subsets are disjoint. The graph $O_n$ has intersection array $(m,m-1,m-1,m-2,m-2, \ldots, \frac{m}{2}+1, \frac{m}{2}+1; 1,1,2,2,3,3, \ldots , \frac{m}{2})$ if $m$ is even, and $(m,m-1,m-1,m-2,m-2, \ldots, \frac{m}{2}+1/2; 1,1,2,2,3,3, \ldots , \frac{m}{2}-1/2, \frac{m}{2}-1/2)$ if $m$ is odd. $O_3$ is commonly known as the Petersen graph, and may be presented in an aesthetically pleasing form:

\vspace{.3in}

\hspace{2.3in}\includegraphics[width=50mm,height=30mm]{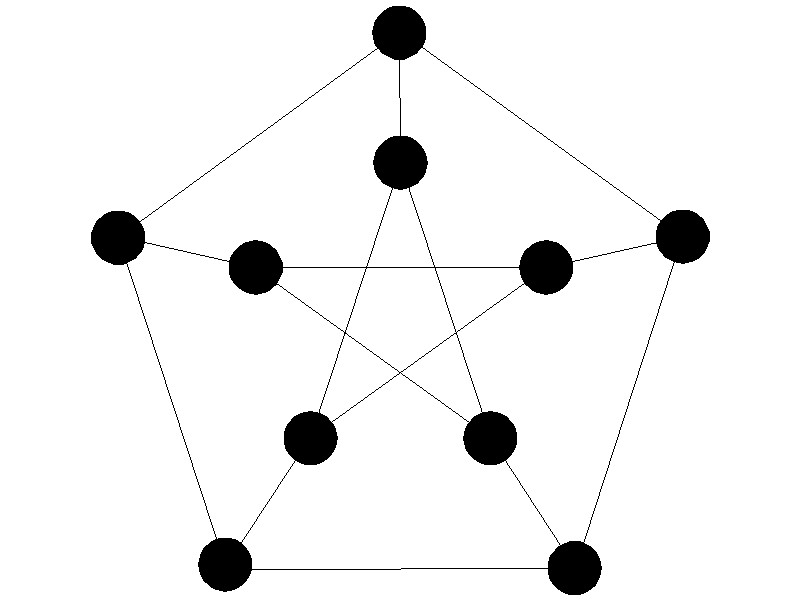}

\vspace{.2in}

(vii) A {\it strongly regular graph} is a regular graph such that there are integers $\la, \mu$ such that any two adjacent vertices have $\la$ common neighbors and any two non-adjacent neighbors have $\mu$ common neighbors. Although the study of these graphs is a field in its own right, the class of strongly regular graphs with $\mu>0$ coincides with the class of distance-regular graphs of diameter 2. An interesting way in which strongly regular graphs differ from distance-regular graphs of higher diameter is that, for strongly regular graphs, $b_1$ may be small relative to the degree $k$. This is important in relation to mixing, hitting, and cover times of random walks; see Theorem \ref{sympdev} and the ensuing discussion for more details.

\vski

(viii) Consider the family of {\it polyhedral graphs}; these are the graphs formed by considering only the vertices and edges of all convex polyhedra in $\RR^m$. Not all polyhedral graphs are distance-regular, but some of the more regular ones are, such as the dodecahedron graph, which is associated with the polyhedron of the same name in $\RR^3$ formed by attaching 12 regular pentagons. The intersection array of the dodecahedron is $(3,2,1,1,1;1,1,1,2,3)$. The following is a presentation of this graph.

\vspace{.3in}

\hspace{2.3in}\includegraphics[width=50mm,height=30mm]{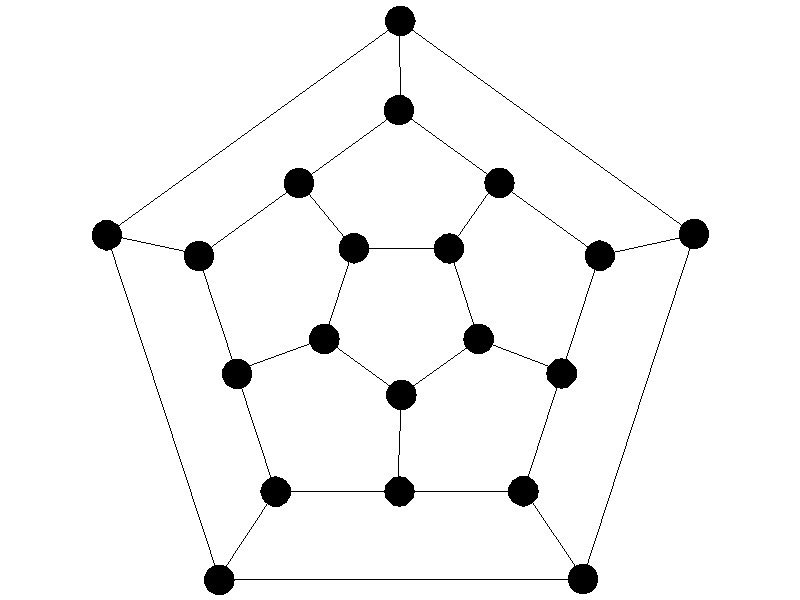}

\vspace{.2in}

(ix) The {\it Biggs-Smith} graph is a degree 3 distance-regular graph with 102 vertices, 153 edges, and intersection array $(3,2,2,2,1,1,1;1,1,1,1,1,1,3)$. It's visualization is difficult, but it is important in the field due to its unusual properties (in \cite{biggsp} it is described as "very exceptional"). For an example of its prominent position within the class of distance-regular graphs see Theorem \ref{bigguy} below.



\vski


There are of course many other examples; a reader interested in seeing an advanced presentation of a larger library of distance-regular graphs is referred to \cite[Ch. 8--14]{drgraphs}.

\section{Electric resistances on distance-regular graphs} \lll{dir}

In this section we discuss the constructions which yield the explicit resistances on distance-regular graphs in terms of the intersection arrays. We will see that there are strong upper and lower bounds on the resistances between points, and that these bounds in turn can yield a great deal of information about simple random walk.

\subsection{The Dirichlet problem with two-point boundaries}
Electric resistance and simple random walk are intimately connected with the concept of harmonic functions on a graph, so we begin there. A function $f$ on $G$ is {\it harmonic} at a point $z\in V$ if $f(z)$ is the average of neighboring values of $f$, that is

\be \label{reff}
\sum_{x\sim z} (f(x)-f(z)) = 0.
\ee

It is easy to see by \rrr{reff} that a function on $G$ harmonic at $z$ satisfies

\be \lll{kyl}
\min_{x\sim z} f(x) \leq f(z) \leq \max_{x\sim z} f(x) ,
\ee

and it follows from this that a function harmonic on all of $G$ must be a constant.
Thus, to obtain interesting functions we must specify a boundary $B$ upon which functions are not required to be harmonic. The {\it Dirichlet problem} is the problem of finding the explicit harmonic function on $G$ with prescribed boundary values. It is easy to see using \rrr{kyl} that a harmonic function must attain its maximum and minimum on $B$, so that a harmonic function with boundary values 0 must be identically 0. By considering the difference in two candidates, this implies that there is always at most one harmonic function with prescribed boundary values; that there always is one can be verified in several ways, but the most relevant for our purposes is, given $f$ defined on $B$, to let $f(x) = E_x[f(X_\tau)]$ on $G \backslash B$, where $\tau$ is the first hitting time of $B$. The Markov property of $X$ shows this function to be harmonic.

\vski

It was Biggs who first noted that the Dirichlet problem can be solved for distance-regular graphs in terms of the intersection array in the case where the boundary contains two points. Suppose $G$ is a distance-regular graph with $n$ vertices and intersection array $(k=b_0,b_1,\ldots ,b_{D-1};$ $c_1,c_2,\ldots,c_D)$. For $0 \leq i \leq D-1$ define the numbers $\phi_i$ recursively by

\bea \label{smile}
&& \pp_0=n-1,
\\ \nn && \pp_i = \frac{c_i\pp_{i-1}-k}{b_i}.
\eea

We will refer to these values as the {\it Biggs potentials}. It can be shown (see \cite{biggsp}) that $\phi_0,\phi_1, \ldots, \phi_{D-1}$ is a strictly decreasing positive sequence. The explicit value of $\phi_i$ is given by the following equation, first stated in \cite{biggsp}:

\be \label{nut}
\phi_i = k\Big(\frac{1}{c_{i+1}} + \frac{b_{i+1}}{c_{i+1}c_{i+2}} + \ldots + \frac{b_{i+1} \ldots b_{D-1}}{c_{i+1} \ldots c_{D}} \Big).
\ee

Fix two adjacent boundary points $u$ and $v$ and let

\bea \label{}
&& K_i^i  = \{x: d(u,x)=i \mbox{ and } d(v,x)=i \},
\\ \nn && K_i^{i+1} = \{x: d(u,x)=i+1 \mbox{ and } d(v,x)=i \},
\\ \nn && K_{i+1}^i = \{x: d(u,x)=i \mbox{ and } d(v,x)=i+1 \}.
\eea
The following fundamental proposition first appeared in \cite{biggsp}.

\begin{proposition} \label{vp}
The function $f$ defined on $G$ by

\bea \label{}
&& f(u) = -f(v) = \pp_0,
\\ \nn && f(x)= 0 \mbox{ for } x \in K_{i}^i,
\\ \nn && f(x)= \pp_i \mbox{ for } x \in K_{i+1}^i,
\\ \nn && f(x)= -\pp_i \mbox{ for } x \in K^{i+1}_i,
\eea

is harmonic on $V-\{u,v\}$.
\end{proposition}

The following diagram, which is of the type popular in the field of distance-regular graphs, may be helpful for understanding the definition of $f$ and the sets $K_j^i$.

\hspace{.8in}\includegraphics[width=130mm,height=110mm]{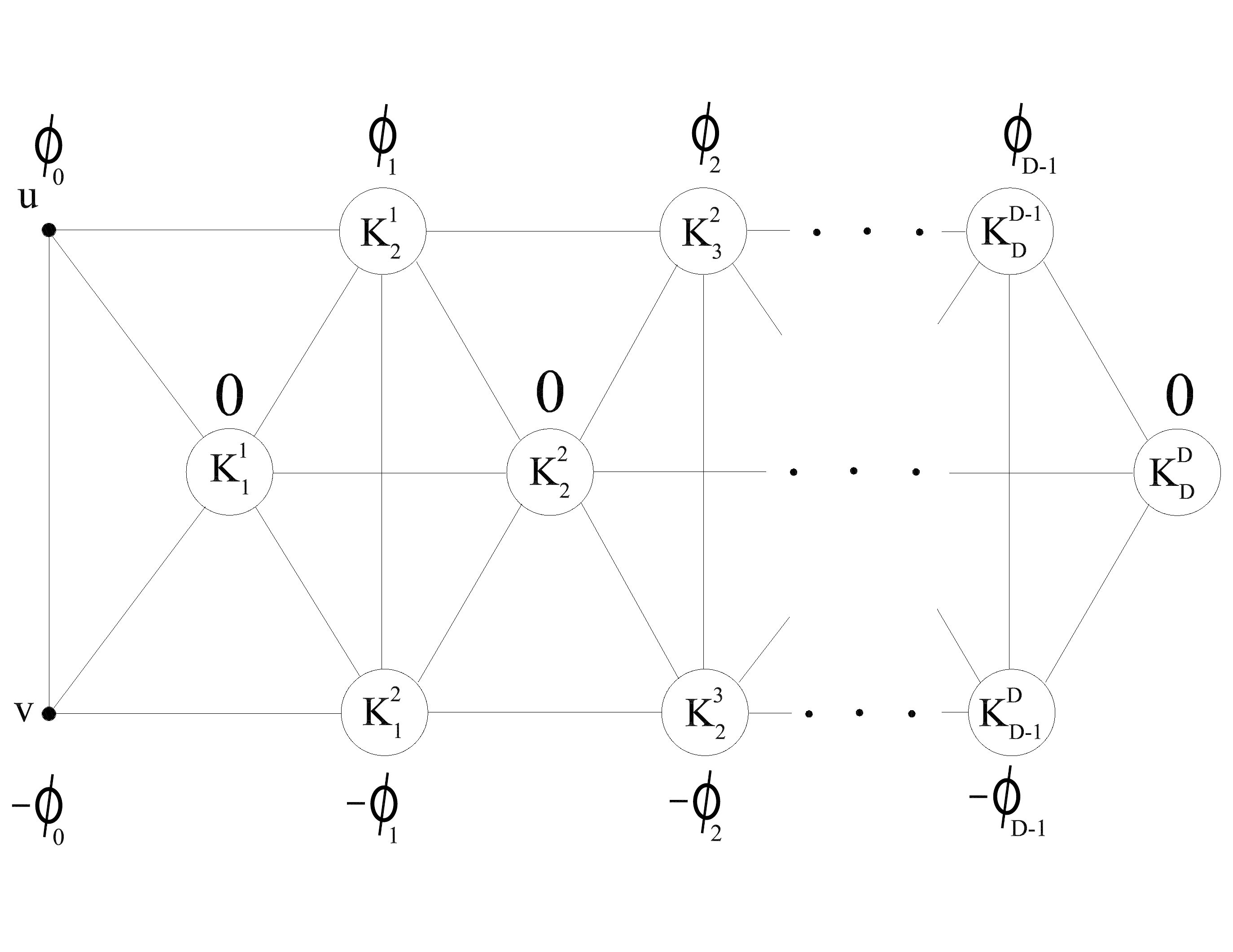}
%

\vspace{-.3in}

Each circle in the diagram represents a set of vertices, with an edge between two circles signifying the fact that each vertex in one circle will in general have neighbors in the other (in particular cases some of these edges may be absent). The proposition can be proved with the aid of the preceding diagram, but as it is subsumed by the more general Proposition \ref{11} below we omit it. By scaling and translation, this function gives the full solution to the Dirichlet problem with two boundary points when the points are adjacent. The $\phi_i$'s can be used to solve the Dirichlet problem when the points are not adjacent as well. This is shown by the ensuing proposition, which is implicit in \cite{biggsp}. In order to simplify notation, for any integer $r$ with $1 \leq r \leq D$ let

\be \lll{}
\Phi_r = \sum_{i=0}^{r-1} \phi_i ,
\ee

and let $\Phi_0=0$.

\begin{proposition} \label{11}
Let $u \neq v$ be vertices of $G$. For any vertex $z$ in $G$ define $f(z) = \Phi_{d(u,z)}, g(z) = \Phi_{d(v,z)}$, and $h(z)=f(z)-g(z)$. Then $h$ is harmonic on $G \backslash \{u,v\}$.
\end{proposition}
The following simple diagram illustrates the definition of $f$, where as before $K_i  = \{z: d(u,z)=i\}$.

\vspace{-.6in}

\hspace{1.1in}\includegraphics[width=110mm,height=90mm]{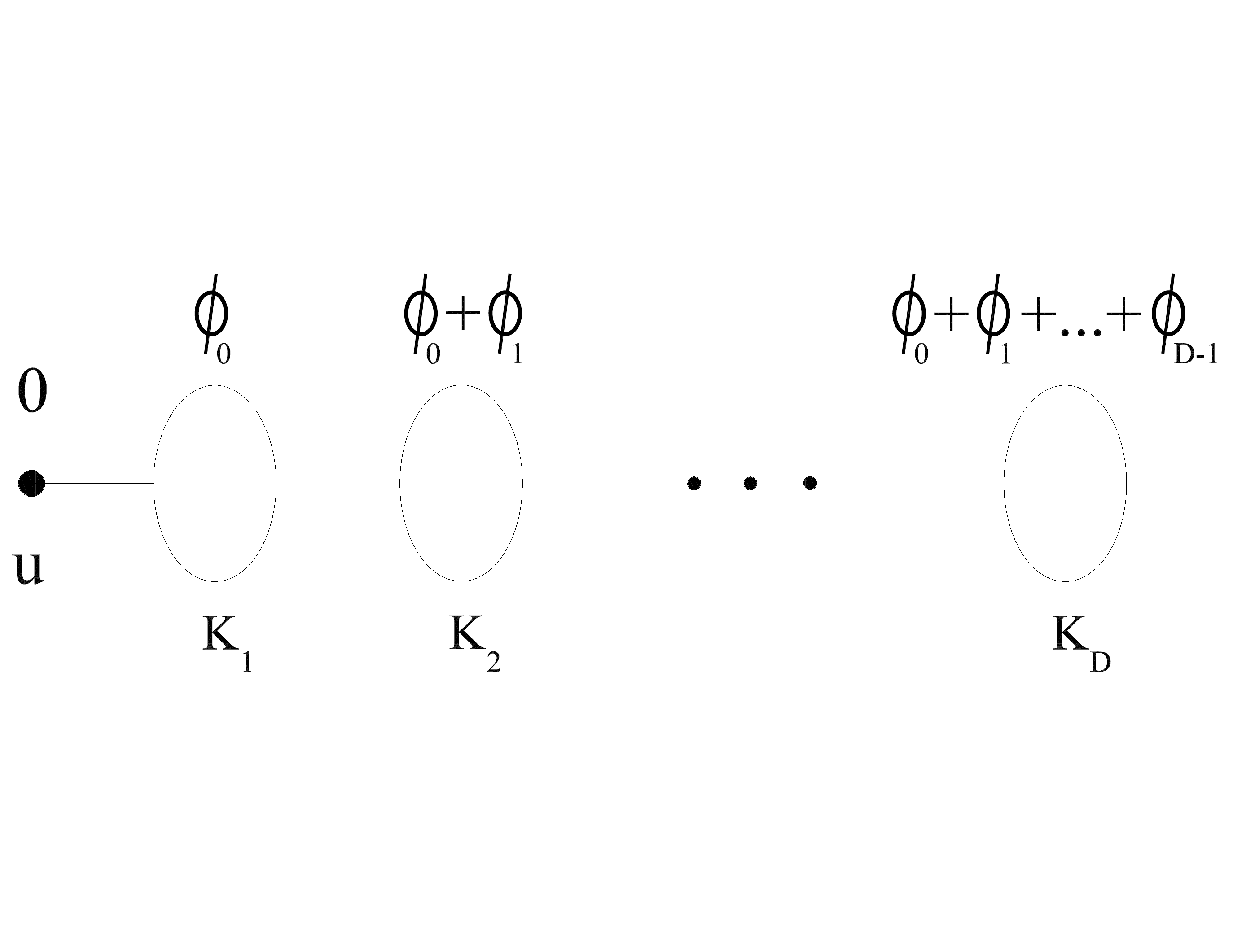}

\vspace{-.8in}


{\bf Proof of Proposition 2:} Let us begin by showing that, for $z \neq u$, we have

\be \lll{truth}
\sum_{x \sim z} f(x) = k f(z) - k.
\ee

Let $r=d(u,z)$. Then

\be \lll{df}
\begin{split}
\sum_{x \sim z} f(x) & = c_r \Phi_{r-1} + a_r \Phi_r + b_r \Phi_{r+1} \\
& = k \Phi_r + b_r \phi_r -c_r \phi_{r-1}.
\end{split}
\ee

But $b_r \phi_r = c_r \phi_{r-1} - k$ by \rrr{smile}, so \rrr{truth} follows. The analogous relation holds for $g$ as well, so that by subtracting these relations we obtain, for $z \neq u,v$,

\be \lll{truth2}
\sum_{x \sim z} h(x) = k h(z).
\ee

Thus, $h$ is harmonic at $z$. \qed

Using this proposition, we will be able to calculate the electric resistance between points, as we now describe. We imagine that $G$ is a circuit where each edge is a wire with resistance 1. We attach a battery of
voltage $V$ to two distinct vertices $u$ and $v$, producing a current through the graph and a voltage $h(x)$ at each vertex. Ohm's Law states that

\be \label{}
V=IR,
\ee

\noi where $V$ represents difference in voltage, $I$ represents current, and $R$ represents resistance. This law applies to each individual edge; that is, the current flowing over each edge is merely the difference in voltage at the endpoints. Ohm's Law can also be applied to the circuit as a whole, which brings us to the concept of the {\it effective resistance} between $u$ and $v$; this is defined as the ratio between the voltage difference $h(v)-h(u)$ and the current flowing from $u$ to $v$, which can be calculated as $\sum_{x \sim u} (h(x)-h(u))$. The other major natural law relevant to us is Kirchhoff's Circuit Law, which states that the net current flowing into any vertex other than $u$ and $v$ is the same as the net current flowing out; that is, $\sum_{x \sim y} (h(x)-h(y)) = 0$ for $y \in G \backslash \{u,v\}$, which means that $h$ is harmonic. The uniqueness of harmonic functions with prescribed boundary values assures us that we may take as our voltage function the $h$ defined in Proposition \ref{11}. It is straightforward to verify that in this case we have $\sum_{x \sim u} (h(x)-h(u)) = nk$, and we arrive at the following proposition, which was stated in \cite{biggsp}


\begin{proposition} \label{kook}
The resistance between two vertices of distance $j$ in $G$ is given by
\be \label{}
\frac{2 \Phi_j}{nk} = \frac{\Phi_j}{m},
\ee

where $m=nk/2$ is the number of edges in $G$.
\end{proposition}

We will see in the next section that this quantity can be given strong bounds for distance-regular graphs, and in the subsequent section that many consequences can be derived for simple random walk.

\subsection{Regularity of electric resistance for distance-regular graphs} \label{treeroll}

Propositions \ref{11} and \ref{kook} above and show us that an understanding of the behavior of the Biggs potentials is crucial for the study of electric resistance on distance-regular graphs. The first major insight in this direction was made by Biggs, who in \cite{biggsp} conjectured the following, which was later proved in \cite{markool}.

\bt[\cite{markool}] \lll{bigguy}
Suppose the degree $k$ of $G$ is greater than $2$. If $d_j$ is the electric resistance between any two vertices of distance $j$ in $G$, then

\be
\max_j d_j = d_D \leq K d_1,
\ee

\noi where $K=1 + \frac{94}{101} \approx 1.931$. Equality
holds only in the case of the Biggs-Smith graph.
\et

{\bf Remark:} This result is false for degree $k=2$, the cycles. To avoid being continuously forced to stipulate "...assuming $k \geq 3$", we will simply assume $k \geq 3$ in what follows, except where stated to the contrary. Furthermore, although this result and further results in this section hold for the complete graph, we will in general assume that the diameter $D>1$, since most of the ensuing results reduce to trivialities for the complete graph.

\vski

Theorem \ref{bigguy} indicates a strong regularity property of distance-regular graphs, as it states that the resistance between any two points is always at most $2 d_1 = \frac{4(n-1)}{nk}< \frac{4}{k}$. In order for this to occur, there must be many possible paths between points which are far from each other. By Proposition \ref{kook} this is equivalent to the following statement:

\be \lll{gh}
\Phi_D=\phi_0+ \ldots + \phi_{D-1} \leq (1 + \frac{94}{101})\phi_0 = (1 + \frac{94}{101})\Phi_1.
\ee

It is not entirely clear what led Biggs to suspect that something like Theorem \ref{bigguy} might be true. In \cite{biggs}, he states only "the conjecture is based on failure to find a counter-example". Nonetheless, \rrr{gh} indicates what turns out to be a key feature of the Biggs potentials, that the sum of the later $\phi_i$'s is dominated by the earlier ones. In particular, we have the following theorems, both proved in \cite{phi1}.

\bt[\cite{markool}] \lll{iii}
For any $m \geq 0$,

\be \label{}
\phi_{m+1} + \ldots + \phi_{D-1} < (3m+3)\phi_m.
\ee
\et

{\bf Remark:} It seems likely that the $3m+3$ can be replaced by a universal constant, but that is just a conjecture at this point.

\bt[\cite{markool}] \lll{p1}
\be
\phi_{2} + \ldots + \phi_{D-1} \leq \phi_1.
\ee

Equality holds only in the case of the dodecahedron.
\et

Theorem \ref{p1} actually implies Theorem \ref{bigguy} rather easily, and is a much stronger statement. It also leads to interesting consequences when coupled with the fact that $\phi_1 \leq \frac{\phi_0}{b_1}$, because strong lower bounds are known on $b_1$ in terms of the degree $k$. In particular, the following is Theorem 16 in \cite{parkool2}.

\begin{theorem}[\cite{parkool2}] \label{sympdev}
Let $\Gamma$ be the subclass of distance-regular graphs consisting of all graphs $G$ such that

\begin{itemize} \label{}

\item $G$ is the line graph of a Morse graph,

\item $G$ is the flag graph of a regular generalized $D$-gon of order $(s, s)$ for some $s$,

\item $G$ is a Taylor graph,

\item $G$ is the Johnson graph $J(7,3)$, or

\item $G$ is the halved 7-cube.
\end{itemize}

Then, if $D >2$ and $G$ is not in $\Gamma$ then $b_1 \geq \frac{k}{2}$. If $D>2$ and $G$ is in $\Gamma$ then $b_1 > \frac{k}{3}$. Finally, if $D=2$ then $b_1 \geq \min (\frac{5k}{16}, \frac{2\sqrt{k}}{1+\sqrt{2}})$.
\end{theorem}

In light of this result, let us define\footnotemark \footnotetext{The cases $k=3,4$ have been completely classified(see \cite[Theorem 7.5.1]{drgraphs} and \cite{brokool}), and the catalog of possible graphs with these degrees is fairly small. The Biggs potentials associated with all distance-regular graphs of degree 3 and 4 can therefore be explicitly calculated, so that the real interest in these results lie in the case $k \geq 5$. This allows us to assume $\min(\frac{4}{k},\frac{94}{101}) = \frac{4}{k}$, simplifying the formula for $C(G,k)$ in the case $D > 2$. Furthermore, for $k \geq 8$ we have $\max(\frac{16}{5k},\frac{1+\sqrt{2}}{2\sqrt{k}}) = \frac{1+\sqrt{2}}{2\sqrt{k}}$.}

\begin{equation} \label{17}
C(G,k) := \begin{cases}\min(\frac{94}{101},\frac{4}{k}) & D>2, G \notin \Gamma,\\
\min(\frac{94}{101},\frac{6}{k}) & D>2 , G \in \Gamma,\\
\max(\frac{16}{5k}, \frac{1+\sqrt{2}}{2\sqrt{k}}) &D=2.
\end{cases}
\end{equation}

\vski

The following corollary of Theorems \ref{bigguy} and \ref{p1} appears in \cite{phi1}.

\begin{corollary} \label{hhh22} In any distance-regular graph of degree $k \geq 3$ with diameter $D > 2$
\be
\sup_{1 \leq i \leq D} \frac{d_i}{d_1} = \frac{d_D}{d_1} \leq 1+C(G,k).
\ee
In particular, $\frac{d_D}{d_1} = 1+O(\frac{1}{k})$ if $D >2$, and $\frac{d_D}{d_1} = 1+O(\frac{1}{\sqrt{k}})$ if $D=2$.
\end{corollary}

This shows that for large $k$, all points become nearly equidistant when measured with respect to the resistance metric. Since resistance is the most important metric in the study of simple random walk, this yields a number of consequences for such walks. We explore this in the next subsection.

\subsection{Consequences for simple random walk} \lll{aml}

In all that follows, we will let $X_t$ be a simple random walk on $G$, and we will use the standard notations $E_x$ and $P_x$ to denote expectation and probability conditioned on $X_0 = x$. For any $x,y \in G$ we let $H_{xy} = E_x[\tau(y)]$, where $\tau(y)$ is the first hitting time of $y$. That is, $H_{xy}$ is the expected number of steps it takes $X_t$ to pass from $x$ to $y$. This is referred to as the {\it hitting time}. The {\it commute time} $C_{xy}$ is the expected
number of
steps necessary for the random walk to travel from $x$ to $y$ and back to $x$, and in the case of distance-regular graphs
is equal to $2H_{xy}$. By Theorem 1 in \cite{comcov},
the expected commute time of a random walk between two points $x$ and $y$ is equal to $2mR_{xy}$, where $R_{xy}$ is the resistance between $x$ and $y$. This quantity can therefore be expressed in terms of the Biggs potentials via Proposition \ref{kook}, specifically we find $H_{xy} = \frac{R_{xy}}{m} = \sum_{i=0}^{d(x,y)-1} \phi_i$, and we may apply Corollary \ref{hhh22} to obtain the following, which shows that the expected hitting time of any two points in a distance-regular graph is essentially $(n-1)$.

\begin{proposition} \label{cook}  Suppose that $G$ is a distance-regular graph. Then for any $x \neq y$ in $G$ we have

\be \lll{freef1}
m \Big( \frac{n-1}{m} \Big) = n-1 \leq H_{xy} \leq m(1+C(G,k)) \Big( \frac{n-1}{m} \Big) = (1+C(G,k))(n-1).
\ee

In particular, $\sup_{x,y\in G} |H_{xy} - (n-1)| = (n-1)O(\frac{1}{k})$ if $D >2$, and $\sup_{x,y\in G} |H_{xy} - (n-1)| = (n-1)O(\frac{1}{\sqrt{k}})$ if $D=2$.
\end{proposition}

\noi The {\it cover time} $Co(G)$ is the expected number of steps that our random walk requires before it has visited every
site on $G$, counting the time $t=0$ as a visit to the initial point. If we let $H^+ = \max_{x,y \in G} H_{xy}$ and $H^- = \min_{x,y \in G, x \neq y} H_{xy}$, then a well-known result due to Matthews \cite{matthews} shows that

\begin{equation} \label{}
H^- \Big( 1 + \frac{1}{2} + \ldots +\frac{1}{n-1}\Big) \leq Co(G) \leq H^+ \Big( 1 + \frac{1}{2} + \ldots +\frac{1}{n-1}\Big).
\end{equation}

It is clear that the precise bounds on $H^+$ and $H^-$ from Proposition \ref{cook} will serve to give precise bounds on $Co(G)$, and, using the simple relations $\log n < 1 + \frac{1}{2} + \ldots +\frac{1}{n-1}$ and $1 + \frac{1}{2} + \ldots +\frac{1}{n-1} < 1+ \log (n-1)$, we have

\begin{proposition} \label{tree}
\begin{equation} \label{kilo}
(n-1) \log n < Co(G) < (n-1)(1+C(G,k))(1+\log (n-1)).
\end{equation}
\end{proposition}

We remark that in \cite{devry} the lower bound in \rrr{kilo} was shown by the same method as here. However, the authors of \cite{devry} state on p. 504 that, in reference to distance-regular graphs, "... it is impossible to derive a general upper bound of the order $n \log n$, since for the cycle graph $Cov_x(G) = n(n-1)/2$" \footnotemark \footnotetext{The notation in the quotation has been changed to that of this paper.}. Nevertheless, Proposition \ref{tree} shows that we do have an upper bound of order $n \log n$ for all but the cycle graphs. We remark further that in \cite{fue}, it was shown that for all graphs, distance-regular or otherwise, we have

\be \label{}
Cov_x(G) \geq n\ln{n}(1+ o(1)).
\ee

Thus, Proposition \ref{tree} shows that it is not possible to find a class of graphs with cover times which grow at an order slower than the distance-regular ones.






%

\vski

Another measure of connectivity of graphs as viewed by random walks was studied in \cite{pal}. There, it was proposed to study $F(u) = \sum_{v \in G} H_{uv}$. It was conjectured that $F(u) \geq (n-1)^2$ for all graphs, independent of $u$, and this was proved for a class of graphs which contains the regular graphs, which of course is a much larger class than the distance-regular graphs. On the other hand, Theorems \ref{bigguy} and \ref{p1} in the form of Proposition \ref{cook} provide a bound in the opposite direction, and shows that $F(u)$ is nearly equal to $(n-1)^2$ for distance-regular graphs of large degree. Combining the results of \cite{pal} with \rrr{freef1} gives

\begin{proposition} \label{aichains} If $G$ is a distance-regular graph, then for every $u$ we have
\begin{equation} \label{}
(n-1)^2 \leq F(u) \leq (1+C(G,k))(n-1)^2.
\end{equation}
In particular, $F(u) = (n-1)^2 (1+O(\frac{1}{k}))$ if $D >2$, and $F(u) = (n-1)^2 (1+O(\frac{1}{\sqrt{k}}))$ if $D=2$.
\end{proposition}

We remark that a number of other consequences concerning the rate of mixing of random walks on distance-regular graphs can be deduced. Let us illustrate this, using the notation of \cite[Ch. 4]{aldfill}. Assume for convenience that $G$ is not bipartite, which implies that the distribution of $X_t$ converges to the stationary distribution, which is the uniform distribution on $G$; if $G$ is bipartite, the ensuing results can be be adjusted by making the walk lazy, as is shown in detail in \cite{levin}. We let $\tau_0 = \frac{1}{n^2} \sum_{u \in G} F(u) = \frac{F(z)}{n}$, where $z$ is any $z \in G$. It is clear that Proposition \ref{aichains} gives the bound $\tau_0 \leq \frac{(1+C(G,k))(n-1)^2}{n}$. For starting point $z \in G$, let $d(t)$ denote the distance of $X_t$ from stationarity, that is,

\begin{equation} \label{reject}
d(t) = \frac{1}{2} \sum_{u \in G} |P_z(X_t = u) - \frac{1}{n}|.
\end{equation}
It is not hard to see that this quantity is independent of $z$ in distance-regular graphs. Let $\tau_1 = \{ t>0 : d(t) < e^{-1} \}$. Let $\tau_2$ be the relaxation time, which is defined by $\tau_2 = 1/(1-\la_2)$, where $\la_2$ is defined to be the second largest eigenvalue of the transition matrix of $X_t$. The importance of this quantity lies in the fact that the $L^2$ distance between the distribution of $X_t$ and the uniform distribution goes to 0 exponentially, and the time constant in the rate of convergence can be given in terms of $\tau_2$. Finally, let

\begin{equation} \label{}
\tau_c = \sup_{A \subseteq G} \frac{|A^c|}{P(X_1 \in A^c|X_0 \mbox{ unif. dist. on }A)}.
\end{equation}
$\tau_c$ is a flow parameter, essentially measuring the difficulty of random walks passing out of $A$ relative to the size of $A^c$. It is shown in \cite[Ch. 4]{aldfill} that a bound on $\tau_0$ gives bounds on the remaining quantities; in particular, we have $\tau_0 \geq \tau_2 \geq \tau_c$, and $66\tau_0 \geq \tau_1$. Thus, the bound on $\tau_0$ contained in Proposition \ref{aichains} yields the following.

\begin{proposition} \label{}
On a distance-regular graph, we have

\begin{equation} \label{}
\begin{gathered}
\tau_1 \leq \frac{66(1+C(G,k))(n-1)^2}{n} \leq 66(1+C(G,k))(n-1), \\
\tau_2 \leq \frac{(1+C(G,k))(n-1)^2}{n} \leq (1+C(G,k))(n-1), \\
\tau_c \leq \frac{(1+C(G,k))(n-1)^2}{n} \leq (1+C(G,k))(n-1).
\end{gathered}
\end{equation}
\end{proposition}

That these bounds grow linearly in $n$ gives quantitative evidence for the statement: "Random walks mix rapidly on distance-regular graphs".

\section{Projected random walks as finite birth-death chains} \lll{proj}

We fix a $v \in G$ and let $Y_t = d(X_t,v)$. $Y_t$ is the random walk projected onto the sets $K_i(v)$, defined in the introduction (the term {\it lumped walk} is also sometimes used). $Y_t$ is a Markov chain on the state space $\{0,1, \ldots, D\}$ with the following transition probabilities:

\be
P(Y_t=j|Y_{t-1} = i) = \left \{ \begin{array}{ll}
\frac{b_i}{k} & \qquad  \mbox{if } j=i+1,    \\
\frac{a_i}{k} & \qquad \mbox{if } j=i, \\
\frac{c_i}{k} & \qquad \mbox{if } j=i-1, \\
0 & \qquad \mbox{if } |i-j|>1. \;
\end{array} \right.
\ee
$Y_t$ is what is known as a {\it birth-death chain} with finite state space. These processes are well understood, which allows us to derive a number of consequences for simple random walks on distance-regular graphs, as the next few subsections show.
%
%
%
%
%
%
%

\subsection{Moments of hitting times and related quantities} \lll{moments}

It is clear that, if we let $\tau_u$ be the first time $X_t$ hits $v$ starting from $u$ and $\und \tau_i$ the first time $Y_t$ hits $0$ starting from $i=d(u,v)$, then $\tau_u = \underline \tau_i$. This allows us to calculate the moments of $\tau_u$, since from the Markov property applied to $Y_t$ we have

\begin{equation} \label{teo}
E[\und \tau_i^q] = \frac{c_i}{k} E[(1+\und \tau_{i-1})^q] + \frac{a_i}{k} E[(1+ \und \tau_{i})^q] + \frac{b_i}{k} E[(1+\und \tau_{i+1})^q].
\end{equation}

This idea was developed by van Slijpe in \cite{veeslee} in order to calculate the variance and first two moments of $\tau_u$. We will see that, using the results on Biggs potentials from Section \ref{treeroll}, we will be able to show that the variance of $\tau_u$ is approximately $(n-1)^2 -(n-1)$. Solving the recurrence relation \rrr{teo} leads to the following theorem. Recall that $k$ is the degree of the graph, and $k_j=|K_j|$, the number of vertices of distance $j$ from any point; we also let $e_j$ denote the number of edges connecting a point in $K_{j-1}$ with a point in $K_j$, which is $e_j = k_j b_{j-1} = k_jc_j$.

\begin{theorem}[\cite{veeslee}] \label{gcrazy}
Suppose $i=d(u,v)$. Then

\begin{equation} \label{}
\begin{gathered}
H_i := E_v[\tau_u] = E[\und \tau_i] = k\sum_{j=1}^{i} \frac{1}{e_j} \sum_{r=j}^{D} k_r, \\
E_v[\tau_u^2] = E[\und \tau_i^2] = -H_i + 2k\sum_{j=1}^{i} \frac{1}{e_j} \sum_{r=j}^{D} k_r H_r.
\end{gathered}
\end{equation}
\end{theorem}

We can cast these expressions in terms of the Biggs potentials by using the following identity, which was proved in Lemma 1 of \cite{markool}:

\begin{equation} \label{}
\phi_{i-1} = \frac{k}{e_i} \sum_{j=i}^{D} k_j.
\end{equation}

We see that $H_i = \sum_{j=0}^{i-1} \phi_j$, which was derived previously in Section \ref{aml}, and

\begin{equation} \label{}
E[\und \tau_i^2] = - \sum_{m=0}^{i-1} \phi_m + 2k\sum_{j=1}^{i} \frac{1}{e_j} \sum_{r=j}^{D} k_r \sum_{s=0}^{r-1} \phi_s.
\end{equation}

Recall now from Proposition \ref{cook} that $(n-1) \leq H_i = \sum_{m=0}^{i-1} \phi_m \leq (1+C(G,k))(n-1)$ for all $i>0$ so that we may bound

\begin{equation} \label{}
\begin{split}
E[\und \tau_i^2] & \leq -(n-1) + 2(1+C(G,k))(n-1) \sum_{j=1}^{i} \frac{k}{e_j} \sum_{r=j}^{D} k_r \\
& = -(n-1) + 2(1+C(G,k))(n-1) \sum_{j=0}^{i-1} \phi_j \\
& \leq  2(1+C(G,k))^2(n-1)^2 - (n-1)
\end{split}
\end{equation}

In a similar fashion, we can give lower bounds for $E[\und \tau_i^2]$, and then similarly give upper and lower bounds for the variance $Var_v(\tau_u) = E_v[\tau_u^2]- H_{vu}^2$. We arrive at the following proposition, which essentially states that $Var_v(\tau_u) \approx (n-1)^2 -(n-1)$.

\begin{proposition} \label{trouble}
Suppose $i=d(u,v)>0$. Then

\begin{equation} \label{}
\begin{gathered}
E_v[\tau_u^2]=E[\und \tau_i^2] \leq 2(1+C(G,k))^2(n-1)^2 - (n-1), \\
E_v[\tau_u^2]=E[\und \tau_i^2] \geq 2(n-1)^2 - (1+C(G,k))(n-1), \\
Var_v(\tau_u)  = E[\und \tau_i^2]- H_i^2 \leq (1+4 C(G,k) + 2C(G,k)^2)(n-1)^2 - (n-1), \\
Var_v(\tau_u) = E[\und \tau_i^2]- H_i^2 \geq (1-2C(G,k)+C(G,k)^2)(n-1)^2 - (1+C(G,k))(n-1). \\
\end{gathered}
\end{equation}
\end{proposition}

van Slijpe also proved a number of interesting identities related to the covering of the graph by a random walk traveling from $v$ to $u$. We will see, again, that our previous work on the Biggs potentials allows us to give strong bounds on these quantities. Let $\tau_u^+$ be the first time $t > 0$ such that $X_t = u$. Then, for $w \neq u$, define

\begin{equation} \label{}
\begin{gathered}
V_{vwu} = \{ X_t = w \mbox{ for some $0 < t \leq \tau_u^+$} \}, \\
M_{vu} = \mbox{\# of different vertices visited at times between 0 and $\tau_u^+$ on a walk started at $v$}, \\
N_{vwu} = \mbox{\# of visits to $w$ at times between 0 and $\tau_u^+$ on a walk started at $v$}.
\end{gathered}
\end{equation}

We then have the following theorem, where we define $H^+_{vu} = E_x[\tau_u^+]$; note that this differs from $H_{vu}$ only when $v=u$, in which case we have $H^+_{vv}=n$, an identity valid for all regular graphs.

\begin{theorem} [\cite{veeslee}] \label{thor}
\begin{equation} \label{}
\begin{gathered}
P(V_{vwu}) = \frac{H^+_{vu} + H^+_{wu} - H^+_{vw}}{2H^+_{wu}}, \\
M_{vu} = \sum_{w \neq u} P(V_{vwu}) = \sum_{W \neq u}\frac{H^+_{vu} + H^+_{wu} - H^+_{vw}}{2H^+_{wu}}, \\
E[N_{vwu}] = \frac{H^+_{vu} + H^+_{wu} - H^+_{vw}}{n}, \\
Var(N_{vwu}) = \frac{(H^+_{vu} + H^+_{wu} - H^+_{vw})(3H^+_{wu} - H^+_{vu} + H^+_{vw} - n)}{n^2}.
\end{gathered}
\end{equation}
\end{theorem}

{\bf Remark:} These identities hold for more general graphs than distance-regular; for the proof in \cite{veeslee}, it was assumed only that $H^+_{xy}=H^+_{yx}$ for all $x,y \in G$. In fact, the first two identities can be realized as special cases of several identities given for Markov chains in \cite[Ch. 2]{aldfill}.

\vski

Again Proposition \ref{cook} allows us to employ the bounds $(n-1) \leq H_{xy} \leq (1+C(G,k))(n-1)$ for all $x,y \in G$, and it is then straightforward to bound the quantities in Theorem \ref{thor}. We obtain the following proposition, which shows that, regardless of the choices of $v,w,u$, we have $P(V_{vwu}) \approx 1/2, M_{vu} \approx (n-1)/2, E[N_{vwu}] \approx (n-1)/n$, and $Var(N_{vwu}) \approx 2(n-1)^2/n^2$.

\begin{proposition} \label{feeling}
\begin{equation} \label{}
\begin{gathered}
\frac{1}{2} - \frac{C(G,k)}{2} \leq P(V_{vwu}) \leq \frac{1}{2} + \frac{C(G,k)}{2}, \\
(n-1)\Big(\frac{1}{2} - \frac{C(G,k)}{2}\Big) \leq M_{vu} \leq (n-1)\Big(\frac{1}{2} + \frac{C(G,k)}{2}\Big), \\
\frac{(1-C(G,k))(n-1)}{n} \leq E[N_{vwu}] \leq \frac{(1+2C(G,k))(n-1)}{n}, \\
\frac{2(1-C(G,k))^2(n-1)^2}{n^2} \leq Var(N_{vwu}) \leq \frac{2(1+2C(G,k))^2(n-1)^2}{n^2}.
\end{gathered}
\end{equation}
\end{proposition}

\subsection{The generating function of hitting times}

Suppose $d(u,v)=i$, and let $GF_i(s) = E_v[s^{\tau_u}]$ be the generating function of the hitting time between two points of distance $i$. In \cite{devry} a formula for $GF_i$ was given, as we now describe. The {\it intersection matrix} of $G$ is defined to be the tridiagonal $(D+1) \times (D+1)$ matrix given by

\begin{equation} \label{}
P =  \left( \begin{array}{ccccccc}
 0 & 1 &  &  &  &  &  \\
 k & a_1 & c_2 &  &  &  &    \\
 & b_1 & a_2 & . &  &  &  \\
  &   & b_2 & . & . &   &  \\
  &   &   & . & . & . &  \\
  &   &   &  & . & . & c_D \\
  &   &   &  &  & . &  a_D \\ \end{array} \right)
\end{equation}

\vski

where all unspecified entries are 0. The {\it adjacency matrix} $\Gamma$ is the $n \times n$ matrix formed by numbering the vertices of $G$ from 1 to $n$ and then letting $\Gamma_{ij} = 1$ if $i \sim j$, and $\Gamma_{ij}=0$ otherwise. The following are very important facts about the algebraic structure of distance-regular graphs (see \cite[Ch. 21]{biggsbook}).

\begin{theorem} \label{cleo}
\begin{itemize}

\item[(i)] $\Gamma$ has exactly $D+1$ distinct eigenvalues $\la_0, \la_1, \ldots, \la_D$, which are precisely the eigenvalues of $B$.

\item[(ii)] Let $u_i, v_i$ be the eigenvectors of $B$ associated to each eigenvalue $\la_i$. Then, for all $i \in \{0, \dots, D\}$, $u_i(0)$ and $v_i(0)$ are both nonzero, and if we standardize such that $u_i(0)=v_i(0)=1$ then

\begin{equation} \label{}
\begin{gathered}
(u_i,v_j) = 0, \mbox{ for $i \neq j$, and} \\
(u_i,v_i) = \frac{n}{m(\la_i)},
\end{gathered}
\end{equation}
where $m(\la_i)$ is the multiplicity of $\la_i$ as an eigenvector of $\Gamma$.
\end{itemize}
\end{theorem}

With this notation, we have

\begin{theorem}[\cite{devry}] \label{varad}
For $i \in \{1, \dots, D\}$ and $|s| \leq 1$, we have

\begin{equation} \label{}
GF_i(s) = \frac{1+(1-s)\sum_{r=1}^{D} m(\la_r)\frac{v_r(i)}{k_i} \frac{1}{1-s\la_r/k}}{1 + (1-s)\sum_{r=1}^{D} m(\la_r)\frac{1}{1-s\la_r/k}}.
\end{equation}
\end{theorem}

Using this result, the authors of \cite{devry} obtained the following proposition on the limiting distributions of cover times on certain distance-regular graphs, among other things. We let $T_{cov}$ be the first time at which every vertex has been visited by our walk, so that the previously defined $Co(G)$ is equal to $E[T_{cov}]$. The {\it standard Gumbel distribution} is the distribution with cumulative distribution function $e^{-e^{-x}}$.

\begin{theorem}[\cite{devry}] \label{clint}

\begin{itemize} \label{}

\item[(i)] Let $H(m,q)$ be the Hamming graph\footnotemark \footnotetext{\cite{devry} uses the term {\it hypercube} for this graph.} as described in Section \ref{hist}, with $m$ fixed. Then, as $q \lar \ff$, $(T_{cov}-n \log n)/n$ converges in distribution to the standard Gumbel distribution.

\item[(ii)] Let $J(m,q)$ be the Johnson graph\footnotemark \footnotetext{\cite{devry} uses the term {\it binomial coefficient graph} for this graph.} as described in Section \ref{hist}. Then, as $m \lar \ff$, $(T_{cov}-n \log n)/n$ converges in distribution to the standard Gumbel distribution, provided that $q=o(m/\log m)$.

\end{itemize}

\end{theorem}

\subsection{Green's function} \lll{}

Let a subset $E$ of $G$ be fixed and let $u \in G \backslash E$. Let $\GG_{u,E}(x)$ denote the expected number of times that $X_t$ starting at $u$ hits $x$ before reaching $E$. We refer to $\GG_{u,E}(x)$ as {\it Green's function with respect to u and E}. The purpose of this section is to show that $\GG_{u,E}(x)$ can be explicitly calculated when $E$ is taken to be a set of uniform distance from $u$. The following lemma is probably well-known, but we provide the short proof for the reader's convenience.

\begin{lemma} \label{}
Let $G$ be a regular graph of degree $k$, $E$ a subset of $G$, and $u \in G \backslash E$. Then $\GG_{u,E}(x)$ is a harmonic function on $G \backslash \{ E, u\}$.
\end{lemma}

{\bf Proof:} Let $X_n^E$ be the walk $X_n$ stopped upon hitting $E$, and assume $x \in G \backslash \{ E, u\}$. Then

\begin{equation} \label{}
\begin{split}
\GG_{u,E}(x) & = \sum_{n=0}^{\ff} P(X^E_n=x) \\
& = \sum_{n=1}^{\ff} \sum_{y \sim x} \frac{1}{k} P(X^E_{n-1}=y) \\
& = \frac{1}{k} \sum_{y \sim x} \sum_{n=0}^{\ff} P(X^E_{n}=y) = \frac{1}{k} \sum_{y \sim x} \GG_{u,E}(y),
\end{split}
\end{equation}
so that $\GG_{u,E}$ is harmonic at $x$. \qed

This shows that the following theorem gives us not only the correct expression for a Green's function, but also another example of a nonconstant harmonic function on a distance-regular graph in terms of the intersection array.

\bt \lll{treeh}
Suppose $u$ is a point in a distance-regular graph $G$ and let us define $\GG_{u,E}$ with $E=\{x:d(u,x)=\aaa\}$. Then, for fixed $x$, let $r=d(u,x)$. We then have, for $0 \leq r < \aaa$,

\be \label{99}
\GG_{u,E}(x) = \frac{k (1+\sum_{j=r+1}^{\aaa-1}\frac{b_{r+1}\ldots b_j}{c_{r+1}\ldots c_j})}{b_r|K_r|}.
\ee

This is the unique harmonic function on $G \backslash \{ E, u\}$ with boundary values 0 on $E$ and $1+\sum_{j=1}^{\aaa-1}\frac{b_{1}\ldots b_j}{c_{1}\ldots c_j}$ at $u$.
\et

{\bf Remark:} Note that $|K_r| = \frac{b_0 \ldots b_{r-1}}{c_1 \ldots c_r}$ for $r \geq 1$, so that \rrr{99} may be expressed entirely in terms of the intersection array of $G$. Expressing $\GG_{u,E}$ in terms of the Biggs potentials in a useful way, however, seems difficult.

\vski

{\bf Proof:} 
Let $\GG^Y_j$ denote the expected number of times that $Y_t=j$ before hitting $K_\aaa$. It is clear by arguments of symmetry that $\GG_{u,E}(x) = \frac{1}{|K_r|} \GG^Y_r$. Note that $P(Y_s = j \mbox{ for some $s > t$}|Y_t = j) = \frac{c_j+a_j + b_j \DD_{j+1}}{k}$, where $\DD_{j+1}$ is the probability that $Y_n$ started at $K_{j+1}$ hits $K_j$ before hitting $K_\aaa$. Thus, for $0 \leq j < \aaa$, $\GG^Y_j$ satisfies

\be \lll{}
\GG^Y_j = 1 + \Big(\frac{c_j+a_j + b_j \DD_{j+1}}{k}\Big)\GG^Y_j.
\ee

This can be rearranged to

\be \lll{88}
\GG^Y_j = \frac{k}{b_j(1-\DD_{j+1})}.
\ee

We need now only compute $\DD_{j+1}$. It is well known\footnotemark
\footnotetext{There are a number of ways to prove \rrr{true}. One way is to write down the correct recurrence relations and verify that the expression on the right side of \rrr{true} is the minimal solution; see \cite{mchain}. Another is to use the aforementioned relationship between electrical resistance and random walks; see \cite{doysne}. Alternatively, the formula for Green's function for a birth-death chain (and thus Theorem \ref{treeh}) may be deduced from Tanaka's formula for Brownian local time; see \cite{markaus} for details.}
that for $0 \leq j < \aa$,

\be \lll{true}
\DD_{j+1} = P(Y_n=K_j \mbox{ before } Y_n=K_\aaa|Y_0 = K_{j+1}) = \frac{\sum_{j=r+1}^{\aaa-1} \frac{b_1\ldots b_j}{c_1\ldots c_j}}{\sum_{j=r}^{\aaa-1} \frac{b_1\ldots b_j}{c_1\ldots c_j}}.
\ee

We therefore have

\begin{equation} \label{}
\begin{split}
\GG_{u,E}(x) &= \frac{\GG^Y_r}{|K_r|} \\
& = \frac{k}{b_j(1-\DD_{j+1})|K_r|} \\
& = \frac{k}{b_r(1-\frac{\sum_{j=r+1}^{\aaa-1} \frac{b_1\ldots b_j}{c_1\ldots c_j}}{\sum_{j=r}^{\aaa-1} \frac{b_1\ldots b_j}{c_1\ldots c_j}})|K_r|} \\
& = \frac{k (1+\sum_{j=r+1}^{\aaa-1}\frac{b_{r+1}\ldots b_j}{c_{r+1}\ldots c_j})}{b_r|K_r|}.
\end{split}
\end{equation}

\qed

\subsection{The cutoff phenomenon} \lll{cutoff}

It is well known that the distribution of an aperiodic irreducible Markov chain $X_t$ on a finite state space approaches a stationary distribution $\tilde{\pi}$ as $t \lar \ff$. For simple random walk on a regular graph this stationary distribution is the uniform distribution. If we let $\hat\pi(j) = \frac{|K_j(u)|}{n}$, so that $\hat\pi$ is the projection of the uniform distribution $\pi$ from $G$ onto $\{0, \ldots, d\}$, then it is not hard to see that $\hat\pi$ is the stationary distribution for $Y_t$. This observation facilitates the study of mixing times on distance-regular graphs. In \cite{bels}, Belsley used the projected walk to establish a result on mixing for a large class of distance-regular graphs. In order to state his result, we need a few definitions. Given two distributions $\nu, \mu$ on the vertices of $G$, let the {\it total variation distance} between them be defined by

\begin{equation} \label{}
||\mu-\nu||_{var} = \frac{1}{2}\sum_{x \in G} |\nu(x)-\mu(x)|;
\end{equation}
the normalization factor $1/2$ is chosen in part so that $||\mu-\nu||_{var} \leq 1$. It is straightforward to verify that if we let $\mu_t, \hat\mu_t$ be the distributions of $X_t, Y_t$, respectively, then $||\mu_t - \pi||_{var} = ||\hat \mu_t - \hat \pi||_{var}$. It is a standard fact, shown by spectral or other methods, that $d(t)=||\mu_t-\pi||$ approaches 0 with exponential speed as $t \lar \ff$ (see \cite[Thm. 4.9]{levin}. However, the constants provided by the general theory are often not particularly sharp. Furthermore, in many cases $d(t)$ stays large until a particular time, at which point it rapidly converges to 0. This time is known as a {\it cutoff time}. For the precise definition, we follow \cite{levin}. For $\eps \in (0,1)$ we let $t_{mix}^G(\eps)$ denote the smallest $t$ such that $d(t) < \eps$ on $G$ (recall that on a distance-regular graph the quantity $d(t)$ does not depend on the starting point of $X_t$). We will say that a sequence $G_n$ of graphs has a cutoff if, for all $\eps \in (0,1)$,

\begin{equation} \label{}
\lim_{n \lar \ff} \frac{t^{G_n}_{mix}(\eps)}{t^{G_n}_{mix}(1-\eps)} = 1.
\end{equation}
Belsley defined a class of distance-regular graphs $\XX$ called the "$q$-examples" which contains the Grassman graphs, sesquilinear forms graphs, dual polar graphs, and half dual polar graphs, among others. He then proved a more precise version of the following theorem.

\begin{theorem}[\cite{bels}] \label{bels}
If $\{G_n\}$ is a sequence of graphs in $\XX$ such that the diameters $D_n$ of $G_n$ approach $\ff$ as $n \lar \ff$, then $\{G_n\}$ displays the cutoff phenomenon, and for any $\eps \in (0,1)$ we have $t^{G_n}_{mix}(\eps) \sim D_n$.
\end{theorem}
The clever proof made use of the fact that the $q$-examples are all graphs for which the ratio $\frac{|K_D|}{|G|}$ is rather large; that is, for which most points in $G$ are of the maximal possible distance from $u$. This means that once the random walk hits the set $K_{D_n}$ it is nearly uniformly distributed on $G_n$, which allows the application of techniques similar to those applied in the context of strong stationary times (see \cite[Ch. 6]{levin}). This property is not shared by all distance-regular graphs, and in particular it was stated in \cite{bels} that certain families of Johnson and Hamming graphs have cutoffs asymptotic to $D_n \log D_n$, so that no general result encompassing all distance-regular graphs seems possible.

\section{Further topics} \lll{ft}

We finish by presenting several additional observations. We begin by discussing formulas for harmonic functions with boundaries with more than two points. We then show how the results on the Biggs potentials can be used to prove that harmonic functions on distance-regular graphs do not stray far from their average values; the resulting inequalities we have referred to as {\it Harnac inequalities}. It is likely that the results in this section can be extended considerably.

\subsection{Dirichlet problem with boundaries with more than two points} \lll{more2}

Let $h_{u,v}$ denote the harmonic function defined in Theorem \ref{df}. Then $h_{u,v}(v)=-h_{u,v}(u)= \Phi_{d(u,v)}$, and by scaling and translation $h_{u,v}$ gives a complete solution to the Dirichet problem with boundary points $u,v$. In theory, the $h$'s so defined give a general solution to the Dirichlet problem with arbitrary boundary, as the following proposition shows. Let $B=\{v_1, \dots, v_q\}$ be a set of vertices in $G$, thought of as a boundary set.

\begin{proposition} \label{}
The functions $h_{v_1,v_2}(z), h_{v_1,v_3}(z), \ldots, h_{v_1,v_q}(z), 1(z)$ form a basis for the set of all functions harmonic on $G \bbs B$, where $1(z)$ denotes the constant function which is equal to $1$ at all points.
\end{proposition}

\pf \rrr{kyl} implies that harmonic functions attain their maximal and minimal values on the boundary. Thus, if two functions $f,g$ agree on $B$ and are harmonic on $G \bbs B$ then by considering $f-g$ we see that $f=g$. Harmonic functions are therefore determined by their boundary values and, as such, the space of harmonic functions on $G \bbs B$ has dimension $q$. Suppose $a_1, \ldots , a_q$ are constants so that

\be \lll{67}
a_1 1(z) + a_2 h_{v_1,v_2}(z) + \ldots + a_q h_{v_1,v_q}(z) = 0,
\ee

for all $z$. It is not hard to see that

\be \lll{}
\sum_{z\in G} h_{v_1,v_i}(z) = 0,
\ee

so by summing \rrr{67} over all $z$ we conclude $a_1=0$. For any $i$,


\be
\sum_{x \sim v_i} (h_{v_1,v_j}(x) - h_{v_1,v_j}(v_i)) \left \{ \begin{array}{ll}
= 0 & \qquad  \mbox{if } i \neq j   \\
\neq 0 & \qquad \mbox{if } i=j \;
\end{array} \right.
\ee

so summing $\rrr{67}$ over all $x \sim v_i$ and subtracting $k(a_2 h_{v_1,v_2}(v_i) + \ldots + a_n h_{v_1,v_q}(v_i))$ shows that $a_i = 0$. Thus, the functions $1(z), h_{v_1,v_2}(z), \ldots, h_{v_1,v_q}(z)$ are linearly independent, and since this set of functions has the correct number of elements it forms a basis. \qed

This shows that in principle all harmonic functions on a distance-regular graph can be expressed in terms of the intersection array by means of the Biggs potentials $\{\phi_i\}_{i=0}^{D-1}$. Define the function $m_{v_i}$ which is harmonic on $G \bbs \{v_1, \ldots , v_q\}$ with $m_{v_i}(v_i)=1$ and $m_{v_i}(v_j)=0$ for $j \neq i$. The collection of such $m$'s will be referred to as {\it harmonic measure} on $B$, and for any function harmonic on $G \bbs B$ we have

\be \lll{}
f(x) = \sum_{i=1}^{q} f(v_i) m_{v_i}(x).
\ee

The functions $\{m_{v_1}, \ldots , m_{v_q}\}$ therefore form an orthonormal basis with respect to the $L^2$ norm on $B$ for the space of harmonic functions on $G \bbs B$. However, in most general cases the formula becomes unwieldy. We will present a number of special cases in which the formula can be handled easily. Let us first suppose $q=3$. For any two points $x,y$, we will use the notation $\Phi_{xy}$ as a shorthand to denote $\Phi_{d(x,y)}$.

\begin{proposition} \label{asla}
Given boundary points $\{u,v,w\}$, the function $m_u(z)$ defined by

\bea \lll{ress}
&& \nn \frac{(\!\Phi_{vw} \! + \! \Phi_{uv} \! - \! \Phi_{uw}\!) \Phi_{wz} \! + \! (\!\Phi_{vw} \! - \! \Phi_{uv} \! + \! \Phi_{uw}\!) \Phi_{vz} \!+ \! (\!-\Phi_{vw}\!+\!\Phi_{uv}\!+\!\Phi_{uw}\!)\Phi_{vw}\! - \!2 \Phi_{vw}\! \Phi_{uz}}
{2(\Phi_{uw}\Phi_{uv} + \Phi_{uv}\Phi_{vw} + \Phi_{uw}\Phi_{vw})-(\Phi_{uw}^2 + \Phi_{uv}^2 + \Phi_{vw}^2)}
\eea

is harmonic on $G\backslash\{u,v,w\}$ and takes the values $1$ at $u$ and $0$ at $v,w$.
\end{proposition}

{\bf Proof:} The denominator of the defining expression of $m_u(z)$ is constant in $z$, so we need only to check harmonicity in the numerator. The numerator can be rewritten as

\bea \lll{} \nn
\Phi_{vw}(\Phi_{wz}-\Phi_{uz}) + \Phi_{vw}(\Phi_{vz}-\Phi_{uz}) + (\Phi_{uv}-\Phi_{uw})(\Phi_{wz}-\Phi_{vz}) + (-\Phi_{vw}+\Phi_{uv}+\Phi_{uw})\Phi_{vw}.
\eea

By Proposition \ref{11}, the functions $\Phi_{wz}-\Phi_{uz}, \Phi_{vz}-\Phi_{uz},$ and $\Phi_{wz}-\Phi_{vz}$ are harmonic on $G\backslash\{u,v,w\}$, and harmonicity of $m_u(z)$ follows. The values at $u,v,w$ can be verified by direct substitution. \qed

For an arbitrary $q>3$, the formula for harmonic measure becomes quite large and difficult to handle. However, if we introduce some symmetry or structure into the boundary set, in certain cases the formula becomes simple. One example is if the boundary is a clique, or more generally a {\it distance-d-clique}, which is defined to be a subset of $G$ with all elements of distance $d$ from each other.

\begin{proposition} \label{asla2}
Suppose $\{u,v_2,\ldots, v_q\}$ is a distance-d-clique. Then the function

\be \label{sux}
\frac{1}{q} + \frac{\sum_{j=2}^{q}\Phi_{v_jz} - (q-1)\Phi_{uz}}{q \Phi_d}
\ee

is harmonic on $G \backslash \{u,v_2,\ldots, v_q\}$ and takes the values 1 at $u$ and $0$ at $v_2, \ldots, v_q$.
\end{proposition}

{\bf Proof:} Similarly to the previous proposition, harmonicity follows from the harmonicity of $(\Phi_{v_jz} - \Phi_{uz})$, and the correctness of the boundary values follow from direct substitution. \qed


\subsection{Harnack inequalities} \lll{har}

Here we give a few of consequences of Theorem \ref{p1} on the regularity of harmonic functions, which we broadly refer to as Harnack inequalities. The classical Harnack inequality for harmonic functions in the complex plane $\CC$ states that if $h(z)$ is a positive harmonic function on the unit disc $\DDD = \{|z| < 1\}$ then

\begin{equation} \label{}
\frac{1-|z|}{1+|z|} h(0) \leq h(z) \leq \frac{1+|z|}{1-|z|} h(0).
\end{equation}

In other words, away from the boundary of $\DDD$, $h(z)$ is close to its average value, $h(0)$. We will show that a similar type of regularity can be deduced for harmonic functions on a distance-regular graph. To simplify the statement of the results, we will assume for the rest of this section that the diameter $D>2$ and $G$ is not in $\Gamma$; there is no difficulty in adjusting the results below to the remaining cases if desired, using Theorem \ref{sympdev}. In this case, Theorem \ref{sympdev} gives us $\frac{\phi_1}{\phi_0} \leq \frac{2}{k}$. Let us begin with the case where the boundary of a harmonic function consists of two points $u$ and $v$. We suppose that we have a function $h(x)$ on $G$ which is harmonic on $G \backslash \{u,v\}$. As in Section \ref{more2}, we let $m_u(x)$ be the unique function on $G$ which is harmonic on $G \backslash \{u,v\}$ and equal to 1 at $u$ and 0 at $v$. We define $m_v(x)$ similarly, except now we set $m_v(v) = 1, m_v(u)=0$. The results in Section \ref{dir} show that

\begin{equation} \label{}
m_u(z) = \frac{1}{2} + \frac{\Phi_{vz} - \Phi_{uz}}{2\Phi_{uv}}, \qquad m_v(z) = \frac{1}{2} + \frac{\Phi_{uz} - \Phi_{vz}}{2\Phi_{uv}}.
\end{equation}
We can now write

\begin{equation} \label{bush}
h(z) = h(u)m_u(z) + h(v)m_v(z) = \frac{h(u)+h(v)}{2} + \frac{(h(u)-h(v))(\Phi_{vz} - \Phi_{uz})}{2\Phi_{uv}}.
\end{equation}

We now easily obtain

\begin{proposition} \label{freef}
For $z \in G \backslash \{u,v\}$ we have

\begin{equation} \label{}
|h(z) -  \frac{h(u)+h(v)}{2}| \leq \frac{|h(u)-h(v)|\phi_1}{\phi_0} \leq \frac{2|h(u)-h(v)|}{k}.
\end{equation}

If $d(x,u),d(x,v) \geq 2$, then

\begin{equation} \label{}
|h(z) -  \frac{h(u)+h(v)}{2}| \leq \frac{|h(u)-h(v)| \phi_1}{2\phi_0} \leq \frac{|h(u)-h(v)|}{k}.
\end{equation}

\end{proposition}

{\bf Proof:} Let us assume that $d(z,v) \geq d(z,u)$. Combining \rrr{bush}, Theorem \ref{p1}, and the fact that $d(z,u),d(z,v) \geq 1$, we have

\begin{equation} \label{howbe}
\begin{split}
|h(z) -  \frac{h(u)+h(v)}{2}| & = \frac{|h(u)-h(v)|}{2\Phi_{uv}}|\Phi_{vz} - \Phi_{uz}| \\
& = \frac{|h(u)-h(v)|}{2\Phi_{uv}}\sum_{i=d(z,u)}^{d(z,v)-1}\phi_i \\
& \leq \frac{|h(u)-h(v)| }{2 \phi_0} \sum_{i=1}^{D-1}\phi_i\\
& \leq \frac{|h(u)-h(v)| (2\phi_1)}{2 \phi_0},
\end{split}
\end{equation}
with the understanding that an empty sum is 0. This yields the result, in conjunction with the aforementioned fact that $\frac{\phi_1}{\phi_0} \leq \frac{2}{k}$. If we assume that $d(z,u),d(z,v) \geq 2$, then the lowest index in the sum in the third line of \rrr{howbe} can be changed to 2, resulting in a final bound of $\frac{|h(u)-h(v)| \phi_1}{2 \phi_0}$. \qed

This shows that harmonic functions on distance-regular graphs with two boundary points are in general quite close to their average value for all their non-boundary points. We can relate this to Proposition \ref{feeling} above by choosing the boundary values $h(u)=1, h(v)=0$, which it can be shown gives $h(z) = P_z(\tau_u < \tau_v)$. Proposition \ref{freef} then shows that $|P_z(\tau_u < \tau_v) - \frac{1}{2}| \leq \frac{2}{k}$ for any disjoint triple $z,u,v$. That is, if our goal is to have a random walk hit $u$ before $v$, then for large $k$ it almost doesn't matter where we start the walk. This is equivalent to the first inequality in Proposition \ref{feeling}.

\vski

Doubtless similar, if more complicated, results can be obtained for larger boundaries, though in the general case we have not investigated this due to the complexity of the formulas and a wont of applications. However, we may derive a simple result for distance-d-cliques, as follows.

\begin{proposition} \label{}
Let $\{v_1, v_2, \ldots, v_q\}$ be a distance-d-clique. Suppose that $h(z)$ is a function on $G$ such that $h$ is harmonic on $G \backslash \{v_1, v_2, \ldots, v_q\}$. Then, for any $z \in G \backslash \{v_1, v_2, \ldots, v_q\}$, we have

\begin{equation} \label{}
|h(z) - \frac{1}{q} \sum_{j=1}^{q}h(v_j)| \leq \frac{4}{k}\Big(\frac{q-1}{q}\Big)\sum_{j=1}^{q}|h(v_j)|.
\end{equation}

If $d(z,v_j) \geq 2$ for all $j$, then we have

\begin{equation} \label{}
|h(z) - \frac{1}{q} \sum_{j=1}^{q}h(v_j)| \leq \frac{2}{k}\Big(\frac{q-1}{q}\Big)\sum_{j=1}^{q}|h(v_j)|.
\end{equation}
\end{proposition}

{\bf Proof:} Recall that $h(z) = \sum_{j=1}^{q} h(v_j) m_{v_j}(z)$, where $m_{v_j}$ is defined as in Proposition \ref{asla2}. Arguing as in Proposition \ref{freef}, we can write

\begin{equation} \label{}
\begin{split}
|m_{v_j}(z) - \frac{1}{q}| & = \Big|\frac{\sum_{1 \leq i \leq q, i \neq j} (\Phi_{v_iz} - \Phi_{v_jz})}{q\Phi_d}\Big| \\
& \leq \frac{(q-1)2 \phi_1}{q \phi_0} \leq \frac{4}{k}\Big(\frac{q-1}{q}\Big).
\end{split}
\end{equation}
This leads to

\begin{equation} \label{}
|h(z) - \frac{1}{q} \sum_{j=1}^{q}h(v_j)| = \Big|\sum_{j=1}^{q} (m_{v_j}(z)- \frac{1}{q})h(v_j) \Big| \leq \frac{4}{k}\Big(\frac{q-1}{q}\Big) \sum_{j=1}^{q}|h(v_j)|.
\end{equation}
For the same reason as in Proposition \ref{freef}, assuming that $d(z,v_j) \geq 2$ for all $j$ allows us to halve the bound in this proposition. \qed

\section{Open questions and concluding remarks}

The following is a list of questions which strike the author as natural and seem to be completely open. Doubtless there are many other worthy questions which may be asked as well.

\begin{itemize}

\item Can the $(3m+3)$ in Theorem \ref{iii} be replaced by the universal constant 2? Would useful consequences follow from this?

\item Can analogs of Theorems \ref{gcrazy} and \ref{thor} and Propositions \ref{trouble} and \ref{feeling} be found for higher moments, or for general moments?

\item Can useful bounds analogous to Theorems \ref{bigguy} and \ref{p1} be found which relate to the expression for the generating function of hitting times given in Theorem \ref{varad}? This would require a deeper understanding of the eigenvalues of $B$ than seems to be currently available.

\item Can large deviations results be found for hitting times, cover times, etc.?

\item Can an analog of Theorem \ref{clint} be proved for all distance-regular graphs?

\item Can bounds analogous to Theorems \ref{bigguy} and \ref{p1} be found which relate to the expression for the Green's function given in Theorem \ref{treeh}?

\item Can the results stated in Section \ref{cutoff} be extended to a larger subclass of distance-regular graphs? Can anything be said which holds for all distance-regular graphs?

\item In reference to Section \ref{ft}, can a formula be found for harmonic measure on arbitrary boundaries, and can Harnack inequalities then be given for functions with these larger boundaries?

\item The quantities considered in Sections \ref{aml} and \ref{moments} are all easily computed for the complete graph; indeed, most of the results here give the correct values if $C(G,k)$ is taken to be 0. In this sense, simple random walk behaves similarly on distance-regular graphs as it does on the complete graph. Can this statement be made more precise? Can distance-regular graphs be considered "almost complete graphs" in the field of probability? For instance, a large number of other random processes have been studied on the complete graph, particularly those which arise in statistical mechanics such as the Potts, Ising, and Heisenberg models. Can known results for these models on the complete graph be extended to distance-regular graphs?

\end{itemize}

It is hoped that distance-regular graphs have been revealed as a natural and interesting setting to study simple random walk, and that other researchers will be motivated to do work in this area as well.

\section{Acknowledgements}

I am very grateful to Jacobus Koolen for introducing me to distance-regular graphs, and for many helpful conversations. I'd also like to thank Tim Garoni and Andrea Collevecchio for useful conversations. I am grateful for support from Priority Research Centers Program through the National Research Foundation of Korea (NRF) funded by the Ministry of Education, Science and Technology (Grant \#2009-0094070), as well as Australian Research Council Grant DP0988483.
\bibliographystyle{alpha}
\bibliography{graph}

\end{document}